\documentclass[12pt]{amsart}
\usepackage[ngerman,english]{babel}
\usepackage{amssymb,latexsym}
\usepackage{amsthm}
\usepackage{implicit_functions}
\begin{document}
\title[Note on expanding implicit functions]{Note on expanding implicit functions\\into formal power series by means of\\multivariable Stirling polynomials}
\author{Alfred Schreiber}
\address{Department of Mathematics\\
         University of Flensburg\\
				 Auf dem Campus 1\\
				 24943 Flensburg, Germany}
\email{agt.schreiber@proton.me}
\urladdr{http://www.alfred-schreiber.de}
\begin{abstract}
Starting from the representation of a function $f(x,y)$ as a formal power series with Taylor coefficients $f_{m,n}$, we establish a formal series for the implicit function $y=y(x)$ such that $f(x,y)=0$ and the coefficients of the series for $y$ depend exclusively on the $f_{m,n}$. The solution to this problem provided here relies on using partial Bell polynomials and their orthogonal companions.
\end{abstract}
\maketitle
\section{Introduction}
The problem of computing the higher derivatives of a function $y=y(x)$, which is implicitly given by an equation $f(x,y)=0$, has been discussed several times already in the mathematical literature of the 19th and 20th centuries. L. Comtet has listed some of these papers in the bibliography of his famous monograph \cite{comt1974}. His own contribution to the problem can be found in \cite{comt1968,comt1974,cofi1974}. Recently, the problem has attracted renewed attention, especially with regard to some of its combinatorial aspects. The results in \cite{cofi1974} have been subjected to careful analysis by Wilde \cite{wild2008}, who also gives new proofs. Zemel \cite{zeme2019} provides an in-depth combinatorial interpretation for those binomial building blocks that appear in the closed formula he proved for the higher derivatives of $y$.

\section{Preliminaries}
The procedure described in the following for calculating the higher derivatives of an implicit function starts from the problem as formulated by Comtet in \cite[p.\,152--153]{comt1974}. There, for a function $f$ given as a formal power series
\[
	f(x,y)=\sum_{m,n\geq0}f_{m,n}\frac{x^m y^n}{m!n!}
\]
(with coefficients $f_{m,n}$ from a fixed commutative field of characteristic zero) Comtet poses the somewhat modified (but equivalent) task of finding a formal power series $y=y(x)=\sum_{n\geq1}y_{n}\frac{x^n}{n!}$ such that \mbox{$f(x,y)=0$.} From this one gets a representation of the $k$-th derivatives $D^{k}(y)$ $(k=1,2,3,\ldots)$ as 
\[
	D^k(y)=y_k+\sum_{n\geq1}y_{k+n}\frac{x^n}{n!}.
\]
In order to be able to compute $y_n=D^n(y)(0)$, we assume $f_{0,0}=0$ and $f_{0,1}\neq0$. Then, by writing
\begin{equation*}
	f(x,y)=\sum_{n\geq0}\varphi_{n}\frac{y^n}{n!},
\end{equation*}
where $\varphi_n=\varphi_{n}(x)=\sum_{m\geq0}f_{m,n}\frac{x^m}{m!}$, we see that $f(x,y)=0$ is equivalent to
\begin{equation}\label{eq1}
	g(y):=\sum_{n\geq1}\varphi_{n}\frac{y^n}{n!}=-\varphi_0.
\end{equation}
The formal power series $g$ is invertible (with respect to $\comp$), since $g(0)=0$ and by assumption $D(g)(0)=\varphi_1=f_{0,1}+x\cdot\sum(\ldots)\neq0$. Let $\inv{g}$ denote the (unique) inverse of $g$. Then, the implicit function $y$ is obtained from \eqref{eq1} in the form 
\begin{equation}\label{eq2}
	y=y(x)=\inv{g}(-\varphi_0(x)).
\end{equation}
Comtet \cite{comt1974} evaluates this expression using the Lagrange inversion formula and determines the coefficients $y_n$ by collecting the terms in $x^n/n!$ that occur in this process. But only in principle! In fact, only some few \textit{ad hoc} calculations are performed that yield explicit representations for $y_1,y_2,y_3,y_4$ (see the table on p.\,153). Of course, this does not tell us what the general coefficient $y_n$ actually looks like. 

In the following, we show how the concepts developed in \cite{schr2015,schr2021a, schr2021b} provide a complete insight into the structure of $y_n$ solely as a function of the coefficients of $f(x,y)$. This is done in two reduction steps.

\section{The first reduction step}
In the first step, we determine the $n$th Taylor coefficient $\inv{g}_n$ of $\inv{g}$: 
\begin{equation}\label{eq3}
	\inv{g}_n=\inv{g}_n(x):=\left[\frac{y^n}{n!}\right]\inv{g}(y)=A_{n,1}(\varphi_1(x),\ldots,\varphi_n(x)),
\end{equation}
where $A_{n,1}$ is the first member of the double-indexed family $A_{n,k}$ of multivariable Stirling polynomials in the indeterminates $X_{1}^{-1},X_2,\ldots,X_{n-k+1}$ with $0\leq k\leq n$; see \cite[Eq.\,(7.2)]{schr2015}. A fundamental (and even characteristic) property of these polynomials is their inverse relationship to the partial Bell polynomials $B_{n,k}$ \cite[Thm.\,5.1]{schr2015}, which states that $\sum_{j=k}^{n}A_{n,j}B_{j,k}=\kronecker{n}{k}$, where $\kronecker{n}{n}=1$, $\kronecker{n}{k}=0$, if $n\neq k$ (Kronecker's symbol). For further information, the reader is referred to the monograph \cite{schr2021b}.
\begin{rem}
In \cite{schr2015,schr2021a,schr2021b}, the collective term `Stirling polynomials of the first and second kind' (in several indeterminates) was proposed for $A_{n,k}$ and $B_{n,k}$ because the associated coefficient sums $A_{n,k}(1,\ldots,1)$ and $B_{n,k}(1,\ldots,1)$ turn out to be just the signed Stirling numbers of the first and the Stirling numbers of the second kind, respectively.
\end{rem}
For our purposes we need the following explicit representation of $A_{n,1}$ as a linear combination of monomial terms \cite[Cor.\,7.2]{schr2015}:
\begin{equation}\label{eq4}
	A_{n,1}\mspace{-2mu}=\mspace{-2mu}X_1^{-(2n-1)}\mspace{-20mu}\sum_{\ptsind{2n-2}{n-1}}\mspace{-2mu}\frac{(-1)^{n-1-r_1}(2n-2-r_1)!}{r_2!\dotsm r_n!(2!)^{r_2}\dotsm(n!)^{r_n}}X_1^{r_1}X_2^{r_2}\dotsm X_n^{r_n}.
\end{equation}
The sum has to be taken over the set $\ptsind{2n-2}{n-1}$ of all partitions of $2n-2$ elements into $n-1$ non-empty blocks, that is, of all sequences $r_1,r_2,\ldots,r_{n}$ of non-negative integers such that $r_1+r_2+\cdots+r_n=n-1$ and $r_1+2r_2+\cdots+nr_n=2n-2$.

From equations \eqref{eq2} and \eqref{eq3} we now get 
\begin{align}
	\label{eq5}
	y(x)&=\sum_{k\geq1}\inv{g}_{k}\frac{(-\varphi_0(x))^k}{k!}\\
	&=\sum_{k\geq1}(-1)^{k}A_{k,1}(\varphi_1(x),\ldots,\varphi_k(x))\frac{\varphi_0(x)^k}{k!}.\notag
\end{align}
Using a well-known property of the partial Bell polynomials (see, for instance, \cite[p.\,133]{comt1974}) and observing that $D^j(\varphi_0)(0)=f_{j,0}$ we have
\begin{equation}\label{eq6}
	\frac{\varphi_0(x)^k}{k!}=\sum_{n\geq k}B_{n,k}(f_{1,0},\ldots,f_{n-k+1,0})\frac{x^n}{n!},
\end{equation}
and thus from \eqref{eq3} and \eqref{eq5}
\begin{equation}\label{eq7}
		y(x)=\sum_{k\geq1}\sum_{n\geq k}(-1)^{k}\inv{g}_{k}(x)B_{n,k}(f_{1,0},\ldots,f_{n-k+1,0})\frac{x^n}{n!},
\end{equation}
where $\inv{g}_k(x)=A_{k,1}(\varphi_1(x),\ldots,\varphi_k(x))$ is well-defined as a formal power series because of $\varphi_1\neq0$. Of course, the term $\inv{g}_k(x)$ hides most of the remaining complexity, which is why we do the following power series `ansatz' in a purely formal way for now:
\[
	A_{k,1}(\varphi_1(x),\ldots,\varphi_k(x))=\sum_{j\geq0}a_{k,j}\frac{x^j}{j!}.
\]
With this we obtain from \eqref{eq7}
\begin{align*}
	y(x)&=\sum_{n,j\geq0}\sum_{k\geq1}(-1)^{k}a_{k,j}B_{n,k}(f_{1,0},\ldots,f_{n-k+1,0})\frac{x^{n+j}}{n!j!}\\
	    &=\sum_{m\geq n\geq 0}\binom{m}{n}\left\{\sum_{k\geq1}(-1)^{k}a_{k,m-n}B_{n,k}(f_{1,0},\ldots,f_{n-k+1,0})\right\}\frac{x^m}{m!}. 	
\end{align*}
Since the coefficient of $x^m/m!$ is nonzero if and only if $m\geq n\geq k$, we obtain the following
\begin{prop}\label{eq8}
\begin{equation}
		y_m=\sum_{n=1}^{m}\binom{m}{n}\left\{\sum_{k=1}^n(-1)^{k}a_{k,m-n}B_{n,k}(f_{1,0},\ldots,f_{n-k+1,0})\right\}.
\end{equation}
\end{prop}
This preliminary result is already suitable to calculate the first coefficients.

\begin{exms}
Let us consider the cases $m=1$ and $m=2$. --- It follows from Proposition \ref{eq8} 
$y_1=(-1)^{1}a_{1,0}B_{1,1}(f_{1,0})=-a_{1,0}f_{1,0}$. Observing $A_{1,1}=X_{1}^{-1}$ we thus obtain $a_{1,0}=\inv{g}_1(0)=A_{1,1}(\varphi_1)(0)=\varphi_1(0)^{-1}=f_{0,1}^{-1}$ and hence $y_1=-f_{1,0}f_{0,1}^{-1}$ which corresponds to the familiar identity $y'(x)=-f_{x}f_{y}^{-1}$.

Already for $m=2$ the computational effort increases noticeably. We have
\begin{align*}
	y_2&=\binom{2}{1}\sum_{k=1}^1(-1)^{k}a_{k,1}B_{1,k}(f_{1,0},\ldots,f_{2-k,0})\\
	&+\binom{2}{2}\sum_{k=1}^2(-1)^{k}a_{k,0}B_{2,k}(f_{1,0},\ldots,f_{3-k,0})\\
	&=-2a_{1,1}B_{1,1}(f_{1,0})-a_{1,0}B_{2,1}(f_{1,0},f_{2,0})+a_{2,0}B_{2,2}(f_{1,0}).
\intertext{Now recall $B_{2,1}=X_2$, $B_{2,2}=X_1^{2}$, $A_{2,1}=-X_1^{-3}X_2$, and observe that $\inv{g}'_1(x)=-\varphi'_1(x)\varphi_1(x)^{-2}$. This yields
}	
	y_2&=-2\inv{g}'_1(0)f_{1,0}-f_{0,1}^{-1}f_{2,0}+\inv{g}_2(0)f_{1,0}^{2}\\
		 &=2\frac{\varphi'_1(0)}{\varphi_1(0)^2}f_{1,0}-f_{0,1}^{-1}f_{2,0}-\frac{\varphi_2(0)}{\varphi_1(0)^3}f_{1,0}^{2}\\
		&=2f_{0,1}^{-2}f_{1,0}f_{1,1}-f_{0,1}^{-1}f_{2,0}-f_{0,1}^{-3}f_{0,2}f_{1,0}^2,
\end{align*}
which of course also follows immediately from $y''(x)=2f_{y}^{-2}f_{x}f_{xy}-f_{y}^{-1}f_{xx}-f_{y}^{-3}f_{yy}f_{x}^2$ if we take $x=0$.
\end{exms}

\begin{rem}
The number of distinct monomials in $D^n(y)$ grows rapidly; it is 9 for $y_3$, 24 for $y_4$, and 91159 for $y_{15}$. Comtet \cite[p.\,175]{comt1974} established a generating function for this sequence and gave a table with some of its values. See also Comtet/Fiolet \cite{cofi1974} and the correction made by Wilde \cite{wild2008}. 
\end{rem}
\section{The second reduction step}
In the second and final step, we will show how the general Taylor coefficient $a_{k,l}$ of $\inv{g}_k(x)$ which appears in Proposition \ref{eq8} can be represented by a polynomial expression depending exclusively on the $f_{m,n}$. As explained in Section 3, we make use of the Stirling polynomials of the first kind to accomplish the series reversion in question. 

To facilitate the evaluation of higher-order derivatives of function powers, we first introduce a simple but useful
\begin{lem}
Let $h$ be a function given by the power series $h(x)=\sum_{n\geq0}h_{n}\frac{x^n}{n!}$, $h_0\neq0$, and let $r,j\in\integers$ with $j\geq0$. Then the following applies:
\begin{equation}\label{auxiliary}
	D^j(h^r)(0)=\sum_{k=0}^j(r)_{k}h_{0}^{r-k}B_{j,k}(h_1,\ldots,h_{j-k+1}).
\end{equation}
\end{lem}
\begin{proof}
We write $h^r$ as a composite function $\idty^r\comp h$ and obtain the following using Fa\`{a} di Bruno's formula \cite[Thm.\,A]{comt1974}:
\[
	D^j(\idty^r\comp h)=\sum_{k=0}^j(D^k(\idty^r)\comp h)\cdot B_{j,k}(D(h),D^2(h),\ldots,D^{j-k+1}(h)).
\]
Now, $(D^k(\idty^r)\comp h)(0)=(r)_{k}h_{0}^{r-k}$, where $(r)_k:=r(r-1)\cdots(r-k+1)$ for $k\geq1$ and $(r)_0:=1$ denotes the falling factorial. Taking into account $D^i(h)(0)=h_i$ ($i=1,2,3,\ldots$), this results in the stated formula \eqref{auxiliary}.
\end{proof}
The expression on the right-hand side of \eqref{auxiliary} can be written more concisely as $\widehat{P}_{j,r}(h_0,h_1,\ldots,h_j)$. Here, $\widehat{P}_{j,r}$ denotes the \emph{potential polynomials} introduced by Comtet \cite[p.\,141, Eq. (5f)]{comt1974}; see also \cite[Eq. (3.10)]{schr2021a}. Note that, for negative $r$, $\widehat{P}_{j,r}$ is a Laurent polynomial in the indeterminates $X_0^{-1},X_1,\ldots,X_j$.

We are now in the position to implement the announced second reduction step.
\begin{prop}
Under the assumptions of Proposition \ref{eq8} and for $m\geq n\geq k\geq 1$ we have
\begin{align*}
a_{k,m-n}=&\sum_{\ptsind{2k-2}{k-1}}\Biggr\lbrace\mspace{-2mu}\frac{(-1)^{k-1-r_1}(2k-2-r_1)!}{r_2!\dotsm r_k!(2!)^{r_2}\dotsm(k!)^{r_k}}\quad\times\sum_{j_1+\cdots+j_k=m-n}\Big\lbrace\frac{(m-n)!}{j_1!\cdots j_k!}\\
           &\quad\times \widehat{P}_{j_1,r_1-2k+1}(f_{0,1},f_{1,1},\ldots,f_{j_1,1})
					 \prod_{\nu=2}^{k}\widehat{P}_{j_\nu,r_\nu}(f_{0,\nu},f_{1,\nu},\ldots,f_{j_\nu,\nu})\Big\rbrace\Biggr\rbrace.\notag
\end{align*}
\end{prop}

\begin{proof}
Since the derivative $D^l$ of $l$-th order is a linear operator for every integer $l\geq0$, we obtain from equation \eqref{eq4}:
\begin{align}
	\label{thm2-1}
	a_{k,l}&=D^l(\inv{g}_k)(0)=D^l(A_{k,1}(\varphi_1,\ldots,\varphi_k))(0)\\
	       &=\sum_{\ptsind{2k-2}{k-1}}\mspace{-2mu}\frac{(-1)^{k-1-r_1}(2k-2-r_1)!}{r_2!\dotsm r_k!(2!)^{r_2}\dotsm(k!)^{r_k}}D^l(\varphi_1^{r_1-2k+1}\varphi_2^{r_2}\dotsm \varphi_k^{r_k})(0).\notag 
\end{align}
We evaluate the term $D^l(\ldots)$ by means of the general Leibniz product rule as follows:
\begin{align}
	\label{thm2-2}
	D^l(\varphi_1^{r_1-2k+1}\varphi_2^{r_2}\dotsm \varphi_k^{r_k})=\mspace{-12mu}\sum_{\substack{j_1+j_2+\cdots+j_k=l \\ 		                                        j_1,j_2,\ldots,j_k\geq0}}\,
	\frac{l!}{j_1!j_2!\cdots j_k!}\;D^{j_1}(\varphi_1^{r_1-2k+1})D^{j_2}(\varphi_2^{r_2})\cdots D^{j_k}(\varphi_k^{r_k}).
\end{align}
Therefore, only expressions like $D^{j_{\nu}}(\varphi_{\nu}^{r_{\nu}})(0)$ remain to be reduced. We first apply the auxiliary statement \eqref{auxiliary} to the first factor on the right-hand side of \eqref{thm2-2}:
\begin{align*}
	D^{j_1}(\varphi_{1}^{r_1-2k+1})(0)&=\widehat{P}_{j_1,r_1-2k+1}(D^0(\varphi_1)(0),D^1(\varphi_1)(0),\ldots,D^{j_1}(\varphi_1)(0))\tag{20\,a}\\
	&=\widehat{P}_{j_1,r_1-2k+1}(f_{0,1},f_{1,1},\ldots,f_{j_1,1}).
\intertext{In a completely analogous manner, we obtain for $\nu\geq 2$:}
	D^{j_\nu}(\varphi_{\nu}^{r_\nu})(0)&=
	\widehat{P}_{j_\nu,r_\nu}(f_{0,\nu},f_{1,\nu},\ldots,f_{j_\nu,\nu}).\tag{20\,b}
\end{align*}
Finally, we obtain the asserted explicit formula for the coefficients $a_{k,m-n}$ by putting $l=m-n$ in \eqref{thm2-1} and \eqref{thm2-2} and combining this with (20\,a,b).
\end{proof}
%

\nocite* 
\end{document}